\documentclass[11pt,a4paper,reqno]{amsart}
\usepackage{amsfonts,amsthm, amscd, epsfig, amsmath, amssymb,enumerate}
\numberwithin{equation}{section}
\newtheorem{Th}{Theorem}[section]
\newtheorem{Le}[Th]{Lemma}

\newtheorem{Pro}[Th]{Proposition}

\theoremstyle{definition}

\newtheorem{remark}[Th]{Remark}

\DeclareMathSymbol{\leqslant}{\mathalpha}{AMSa}{"36} 
\DeclareMathSymbol{\geqslant}{\mathalpha}{AMSa}{"3E} 
\DeclareMathSymbol{\eset}{\mathalpha}{AMSb}{"3F}     
\renewcommand{\leq}{\;\leqslant\;}                   
\renewcommand{\geq}{\;\geqslant\;}                   
\newcommand{\suptwo}[2]{\sup_{\substack{#1 \\ #2}}} 
\newcommand{\inftwo}[2]{\inf_{\substack{#1 \\ #2}}} 
\newcommand{\sumtwo}[2]{\sum_{\substack{#1 \\ #2}}} 

\makeatletter
\def\captionfont@{\footnotesize}
\def\captionheadfont@{\scshape}

\long\def\@makecaption#1#2{%
  \vspace{2mm}
  \setbox\@tempboxa\vbox{\color@setgroup
    \advance\hsize-6pc\noindent
    \captionfont@\captionheadfont@#1\@xp\@ifnotempty\@xp
        {\@cdr#2\@nil}{.\captionfont@\upshape\enspace#2}%
    \unskip\kern-6pc\par
    \global\setbox\@ne\lastbox\color@endgroup}%
  \ifhbox\@ne 
    \setbox\@ne\hbox{\unhbox\@ne\unskip\unskip\unpenalty\unkern}%
  \fi
  \ifdim\wd\@tempboxa=\z@ 
    \setbox\@ne\hbox to\columnwidth{\hss\kern-6pc\box\@ne\hss}%
  \else 
    \setbox\@ne\vbox{\unvbox\@tempboxa\parskip\z@skip
        \noindent\unhbox\@ne\advance\hsize-6pc\par}%
\fi
  \ifnum\@tempcnta<64 
    \addvspace\abovecaptionskip
    \moveright 3pc\box\@ne
  \else 
    \moveright 3pc\box\@ne
    \nobreak
    \vskip\belowcaptionskip
  \fi
\relax
}
\makeatother
\def\writefig#1 #2 #3 {\rlap{\kern #1 truecm
\raise #2 truecm \hbox{#3}}}

\newcommand{\cA}{\ensuremath{\mathcal A}}
\newcommand{\cB}{\ensuremath{\mathcal B}}

\newcommand{\cD}{\ensuremath{\mathcal D}}
\newcommand{\cE}{\ensuremath{\mathcal E}}
\newcommand{\cF}{\ensuremath{\mathcal F}}
\newcommand{\cG}{\ensuremath{\mathcal G}}
\newcommand{\cH}{\ensuremath{\mathcal H}}

\newcommand{\cL}{\ensuremath{\mathcal L}}
\newcommand{\cM}{\ensuremath{\mathcal M}}

\newcommand{\cP}{\ensuremath{\mathcal P}}
\newcommand{\cQ}{\ensuremath{\mathcal Q}}
\newcommand{\cR}{\ensuremath{\mathcal R}}
\newcommand{\cS}{\ensuremath{\mathcal S}}

\newcommand{\frH}{\ensuremath{\mathfrak H}}

\newcommand{\bbC}{{\ensuremath{\mathbb C}} }

\newcommand{\bbH}{{\ensuremath{\mathbb H}} }

\newcommand{\bbN}{{\ensuremath{\mathbb N}} }

\newcommand{\bbR}{{\ensuremath{\mathbb R}} }

\newcommand{\bbZ}{{\ensuremath{\mathbb Z}} }

\newcommand{\Om}{\Omega} \newcommand{\om}{\omega}
 \newcommand{\var}{{\rm Var}}
\newcommand{\be}{\begin{equation}}
\newcommand{\bestar}{\begin{equation*}}
\newcommand{\la}{\label} 
\newcommand{\grad}{\nabla} \newcommand{\si}{\sigma}
 
\newcommand{\al}{\alpha} \newcommand{\ket}[1]{\left\vert#1\right\rangle}

\newcommand{\scalar}[2]{\left\langle #1,#2 \right\rangle}

\newcommand{\tg}{\tilde \gamma}

\newcommand{\tfi}{\tilde \varphi}

\let\a=\alpha    \let\d=\delta  \let\e=\varepsilon
 \let\g=\gamma \let\h=\eta      \let\l=\lambda
      \let\o=\omega      
\let\r=\rho      
  
   \let\G=\Gamma  \let\L=\Lambda

\def\\{\hfill\break}

\def\thsp{\thinspace}

\def\tthsp{\kern .083333 em}

\def\?{\mskip -10mu}
\def\indbox#1{\hbox to \parindent{\hfil\ #1\hfil} }

\newfam\msafam
\newfam\msbfam
\newfam\eufmfam
\def\hexnumber#1{%
  \ifcase#1 0\or 1\or 2\or 3\or 4\or 5\or 6\or 7\or 8\or
  9\or A\or B\or C\or D\or E\or F\fi}
\font\tenmsa=msam10
\font\sevenmsa=msam7
\font\fivemsa=msam5
\textfont\msafam=\tenmsa
\scriptfont\msafam=\sevenmsa
\scriptscriptfont\msafam=\fivemsa
\edef\msafamhexnumber{\hexnumber\msafam}%
\mathchardef\restriction"1\msafamhexnumber16
\mathchardef\ssim"0218
\mathchardef\square"0\msafamhexnumber03
\mathchardef\eqd"3\msafamhexnumber2C
\def\QED{\ifhmode\unskip\nobreak\fi\quad
  \ifmmode\square\else$\square$\fi}
\font\tenmsb=msbm10
\font\sevenmsb=msbm7
\font\fivemsb=msbm5
\textfont\msbfam=\tenmsb
\scriptfont\msbfam=\sevenmsb
\scriptscriptfont\msbfam=\fivemsb

\font\teneufm=eufm10
\font\seveneufm=eufm7
\font\fiveeufm=eufm5
\textfont\eufmfam=\teneufm
\scriptfont\eufmfam=\seveneufm
\scriptscriptfont\eufmfam=\fiveeufm

\def\({\left(}
\def\){\right)}
\let\neper=e
\let\ii=i

\def\ie{\hbox{\it i.e.\ }}
\let\id=\identity

\def\nep#1{ \neper^{#1}}

\def\ov#1{{1\over#1}}

\def\tc{\thsp | \thsp}
\let\<=\langle
\let\>=\rangle

\def\gap{\mathop{\rm gap}\nolimits}

\outer\def\nproclaim#1 [#2]#3. #4\par{\medbreak \noindent
   \talato(#2){\bf #1 \Thm[#2]#3.\enspace }%
   {\sl #4\par }\ifdim \lastskip <\medskipamount
   \removelastskip \penalty 55\medskip \fi}
\def\thmm[#1]{#1}
\def\teo[#1]{#1}
\def\sttilde#1{%
\dimen2=\fontdimen5\textfont0
\setbox0=\hbox{$\mathchar"7E$}
\setbox1=\hbox{$\scriptstyle #1$}
\dimen0=\wd0
\dimen1=\wd1
\advance\dimen1 by -\dimen0
\divide\dimen1 by 2
\vbox{\offinterlineskip%
   \moveright\dimen1 \box0 \kern - \dimen2\box1}
}
\begin{document}
\title[Asymmetric diffusions] {Relaxation time of
anisotropic simple exclusion processes and quantum Heisenberg models}
\date{\today} \author[P. Caputo]{Pietro Caputo}
\address{Dip. Matematica, Universita' di Roma Tre, L.go S. Murialdo 1,
00146 Roma, Italy} \email{caputo\@@mat.uniroma3.it}
\author[F. Martinelli]{Fabio Martinelli} \address{Dip. Matematica,
Universita' di Roma Tre, L.go S. Murialdo 1, 00146 Roma, Italy}
\email{martin\@@mat.uniroma3.it}

\vskip 1cm
\begin{abstract}
Motivated by an exact mapping between anisotropic half integer spin
quantum Heisenberg models and asymmetric diffusions on the lattice, we
consider an anisotropic simple exclusion process with $N$ particles in a
rectangle of $\bbZ^2$. Every particle at row $h$ tries to jump to an
arbitrary empty site at row $h\pm 1$ with rate $q^{\pm 1}$, where $q\in
(0,1)$ is a measure of the drift driving the particles towards the
bottom of the rectangle. We prove that the spectral gap
of the generator is uniformly positive in $N$ and in the size of the
rectangle. The proof is inspired by a recent interesting technique
envisaged by E. Carlen, M.C. Carvalho and M. Loss to analyze the Kac
model for the non linear Boltzmann equation.  We then apply the result
to prove precise upper and lower bounds on the energy gap for the
spin--S, ${\rm S}\in \ov2\bbN$, XXZ chain and for the 111 interface
of the spin--S XXZ Heisenberg model, thus generalizing previous
results valid only for spin $\ov2$.

\vskip1cm



\end{abstract}

\maketitle



\section{Introduction}

\la{exclusion}

\noindent
Some years ago it was discovered by Alcaraz \cite{Alcaraz} that a class
of asymmetric reversible simple exclusion processes on $\bbZ^d$ related
to models of diffusion limited chemical reactions, are unitarily
equivalent to certain anisotropic quantum Heisenberg Hamiltonian, known
as XXZ models, that have received in recent years increasing attention
in connection with the analysis of quantum domain walls (see
\cite{Alcarazetal}, \cite{B1}, \cite{S} and references therein). Such an
equivalence implies that the spectrum of (minus) the Markov generator of
the process coincides with the spectrum of the quantum Hamiltonian. In
particular the energy gap above the quantum ground state, a key quantity
in the theory of quantum spin systems, becomes identical to the spectral
gap of the process and a variety of probabilistic techniques come into play
in order to obtain meaningful estimates. Such an observation was exploited
recently in \cite{CapMar} to prove sharp bounds on the energy of low
lying excitations above a ground state describing a 111 interface for a
spin $\ov2$ XXZ model. The extension of such results to higher half
integer spin ${\rm S} \in \ov2 \bbN$ requires additional finer analysis
and led us to consider the following model.

\subsection{Setup}
Given two natural
numbers $L,H$
we consider the rectangle 
\be
\L=\{(i,h)\in \bbZ^2:\; i=1,\dots,L \;{\rm and}\; h=1,\dots,H\}
\la{Lambda}
\end{equation}
For each $i$, $\L_i$ stands for the {\em stick} at $i$ given
by $\L_i = \{(i,h): \; h=1,\dots,H\}$. At each $x\in\L$ we have
a variable $\al_x\in\{0,1\}$: we say that site $x$ is occupied
(by a particle) if $\al_x=1$ and empty otherwise. The set of
configurations $\{0,1\}^\L$ is denoted by $\Om$ and it is naturally
decomposed in single stick configurations: $\al\in\Om$ will
be written often in the form $\al=(\eta_1,\dots,\eta_L)$ with
$\eta_i\in\{0,1\}^H$ denoting the restriction of $\al$ to the
stick $\L_i$.

\smallskip
Given a parameter $q\in(0,1)$ we
define the product probability measure $\mu$
on $\Om$.
\begin{equation}
\mu(f) = \sum_{\al\in\Om}\mu(\al)f(\al)\,,\quad \mu(\al)
=\prod_{i=1}^L \prod_{h=1}^H
\frac{q^{2h\al_{(i,h)}}}{1+q^{2h}}
\la{mula}
\end{equation}
where $f$ is a generic function $f:\Om\to \bbR$.
According to $\mu$ particles prefer to live on the region
of small $h$, i.e.\ the bottom of
the box $\L$ if we interpret $h$ as a vertical coordinate.
We define $n_i$ as the number of particles in the stick $\L_i$:
$n_i(\al)=n_i(\eta_i)=\sum_{h=1}^H\al_{(i,h)}$ and consider the
conditional probability measure
\be
\nu = \nu_{N} = \mu\Big(\,\cdot \tc \sum_{i=1}^L n_i = N\Big)
\la{canon}
\end{equation}
The variance of a function $f$ w.r.t.\ $\nu$ will be written
as usual in one of the following ways
$$
\var(f) = \nu(f,f) = \nu\big((f - \nu(f))^2\big)\,.
$$

\subsection{The process and main result}
The asymmetric diffusion that will be analyzed in the sequel can be
described as follows. Every particle at row $h$ tries to jump to an
arbitrary empty site at row $h+1$ with rate $q$ and to an empty site at
row $h-1$ with rate $1/q$.  The Markov generator is defined by the
operator \be \cL f(\al) = \frac1L \sum_{i=1}^{L}\sum_{j=1}^{L}
\sum_{h=1}^{H-1} c_{(i,h);(j,h+1)}(\al)\grad_{(i,h);(j,h+1)}f(\al)
\la{genera}
\end{equation}
where we use the notation
\be
\grad_{(i,h);(j,h+1)}f(\al) = f(\al^{(i,h);(j,h+1)})-f(\al)\,,\;
\la{grads}
\end{equation}
$\al^{(i,h);(j,h+1)}$ denoting the configuration
in which the values of $\al$ at ${(i,h)}$ and $(j,h+1)$
have been interchanged while the rest is kept unchanged.
The rates $c_{(i,h);(j,h+1)}$ are given by
\be
c_{(i,h);(j,h+1)}(\al)= q^{\al_{(i,h)}-\al_{(j,h+1)}}\,.
\la{rates}
\end{equation}
Simple computations show that $\cL$ is self adjoint in $L^2(\nu)$,
the associated Dirichlet--form being
\begin{gather}
\cD(f,f) =
\nu(f(-\cL)f) = \frac1{L} \sum_{i=1}^L\sum_{j=1}^{L} D_{ij}(f)
\la{dirD0} \\
D_{ij}(f) := \frac12 \sum_{h=1}^{H-1}\nu \big[
c_{(i,h);(j,h+1)}\big(\grad_{(i,h);(j,h+1)}f\big)^2\big]
\nonumber
\end{gather}
Our main result then says that decay to equilibrium for the dynamics
defined by (\ref{genera}) occurs exponentially fast in the $L^2(\nu)$
norm, uniformly in $L,H$ and $N$.
As a corollary we shall obtain
an interesting estimate on the energy gap for a class
of quantum XXZ Hamiltonian, see Theorems \ref{gaps} and \ref{hd} below,
that extends in particular previous results
in \cite{CapMar} and \cite{BKS}.
\\
More precisely let
\be
\gamma(L,H) = \sup_{N} \sup_{f\in L^2(\nu)}\,
\frac{\var(f)}{\cD(f,f)}
\la{gammalh}
\end{equation}
where the number of particles $N$ in $\sup_N$, using the particle--hole
symmetry, is assumed to range from $1$ to $\frac{LH}{2}$.
\begin{Th}
\la{main}
For every $q\in(0,1)$ there exists $C<\infty$ such that
$$\sup_{L,H}\gamma(L,H)\leq C\,.$$
\end{Th}
\begin{remark}
Recently we learned \cite{BBHM} that for $L=1$ and $N= \frac{H}{2}$ the
mixing time (namely the smallest time such that
$\sup_{\a,\tilde\a}\nu\left(|\nep{t\cL}(\a,\cdot)-
\nep{t\cL}(\tilde\a,\cdot)|\right) \le \ov2$) grows like $H$. Remarkably 
in the same setting the logarithmic Sobolev constant grows at least like $H^2$.
\end{remark}
\subsection{Applications}
Some of the applications of Theorem \ref{main} we have in mind,
particularly those to quantum Heisenberg models, are linked to the analysis
of the restriction of the above defined process to the horizontal sums of
the basic variables $\a_{i,h}$ given by
$$
\o_h = \sum_{i=1}^L \al_{(i,h)}\,,\quad\quad h=1,\dots,H
$$
We will show below that the evolution of the new variables
$\{\o_h\}_{h=1}^H$ is still Markovian and that it can be interpreted as
describing the fluctuations of a non--negative profile subject to a
fixed area constraint.
\\
Let $\cP_L$ denote the set of permutations of $\{1,\dots,L\}$.
Given $\pi\in\cP_L$ we write $\a^{\pi,h}$
for the configuration
\be
\a^{\pi,h}_{(i,h')}=\begin{cases}
\a_{(i,h')} & h'\neq h \\
\a_{(\pi(i),h)} & h'=h
\end{cases}\,,\quad\quad i=1,\dots,L
\la{symmpi}
\end{equation}
The subspace $\cS$ of horizontally symmetric functions is defined by
\be
\cS = \{f\in L^2(\nu):\; f(\al)=
f(\al^{\pi,h})\,,\; \forall \pi\in\cP_L,\;\forall h=1,\dots,H\}\,.
\la{symmo}
\end{equation}
and it clearly consists of functions which
only depend on the horizontal sums.
\\
Given $f\in\cS$ we write $\hat f(\o) = f(\a)$. In this way we identify
$\cS$ with the space $L^2(\hat\Om,\hat \nu)$, $\hat \Om = \{0,1,\dots,L+1\}^H$
and $\hat\nu$ the marginal of $\nu$ on horizontal sums $\o=\{\o_h\}$.
The probability $\hat\nu(\o)$ of a single $\o\in\hat\Om$
compatible with the global constraint $\sum_h \o_h = N$
is easily computed to be
\be
\hat\nu(\o) = \frac1{Z}\,\prod_{h=1}^H\binom{L}{\o_h}q^{2h\o_h}
\,\quad Z = \sumtwo{\o\in\hat \Om:}{\sum_h\o_h = N}
\prod_{h=1}^H\binom{L}{\o_h}q^{2h\o_h}
\la{hatnuo}
\end{equation}
Observe that $\cS$ is an invariant subspace for the generator $\cL$,
i.e.\ $\cL\cS\subset \cS$. In fact a simple computation shows that
for every $f\in\cS$, $\cL f(\a)$ can be written as
\begin{gather}
\widehat{\cL} \hat f(\o)= \frac1L
\sum_{h=1}^{H-1}
\Big\{\,w_{+,h}(\o)\big[\hat f(\o^{+,h}) -\hat f(\o)\big]
+ w_{-,h}(\o) \big[\hat f(\o^{-,h}) -\hat f(\o)\big]\,\Big\}
\la{generas}
\\
w_{+,h} :=
q^{-1}\,(L-\o_h)\o_{h+1}\,,\quad
w_{-,h} :=
q\,(L-\o_{h+1})\o_{h}\,,\quad \nonumber \\
\o^{\pm,h}_{h'} :=  \begin{cases}
\o_{h'} & h'\neq h,h+1 \\
\o_{h}\pm 1 & h'=h \\
\o_{h+1} \mp 1 & h'=h+1
\end{cases}
\nonumber
\end{gather}
This defines a Markov generator $\widehat{\cL}$ which is symmetric in
$L^2(\hat\Om,\hat\nu)$.  The corresponding process (the restriction to
$\{\o_h\}$ of the original anisotropic exclusion dynamics) can be
interpreted as describing fluctuations of a non-negative profile
$\o:=\{\o_h\}_{h=1}^H$ subject to a fixed area constraint ($\sum_h \o_h
=$ constant).  In view of the anisotropy the profile is rather sharply
localized: letting $\r=N/L$ we see that $\om_h\approx L$ for heights $h$
below $\r$ and $\o_h\approx 0$ above $\r$ with high probability.  By
Theorem \ref{main} relaxation to equilibrium in $L^2(\hat\Om,\hat\nu)$
is exponentially fast uniformly in $\r$.  \\ In the case $L=2$ the
previous analysis admits another interesting interpretation as a model
for diffusion limited chemical reactions, see \cite{Alcaraz,BG,BL} and
references therein.  Namely describe the state $\o_h=2$ as the presence
at $h$ of a particle of type $A$, $\o_h=0$ as a particle of type $B$ and
$\o_h=1$ as the absence of particles ({\em inert}). If $n_A$, $n_B$
denote the size of the two populations we see that the difference
$n_A-n_B$ is conserved and this system can be studied as a model for
asymmetric diffusion with creation and annihilation of the two
species. Particles of type $A$ have a constant drift towards the bottom
(``small $h$'' region) while particles of type $B$ have the same drift
towards the top (``large $h$'' region). They perform asymmetric simple
exclusion with respect to the {\em inert} sites but when they meet (\ie
when they become nearest neighbors) they can produce the annihilation
reaction $A\,+\,B\,\rightarrow\,$ {\em inert}. The
reverse reaction {\em inert} $\,\rightarrow\,A\,+\,B$
restores steady state fluctuations given by the
canonical measure.
\\
While Theorem \ref{main} implies immediately
$L^2$-exponential ergodicity for the above process, a direct proof
of the result for the two-particle model seemed difficult to us.

\subsection{Some ideas for the proof of Theorem \ref{main}}
We conclude this introductory section with some comments on the main
ideas behind the proof of Theorem \ref{main}. Our main source of
inspiration has been a recent work by E. Carlen, M.C. Carvalho and
M. Loss on the rate of approach to equilibrium for the Kac model of the
non linear Boltzmann equation \cite{CCL}. Like in other approaches to bound the
spectral gap for large reversible Markov chains, the first idea is to recursively
bound $\g(L,H)$ in terms $\g(1,H)$, the latter being finite uniformly in
$H$ because of Theorem 4.3 in \cite{CapMar}.  The starting point, as
e.g. in the martingale approach of H.-T. Yau \cite{Yau}, is a decomposition of
the variance of an arbitrary function $f$ as
\begin{equation*}
\var(f) = \frac1L\sum_{k=1}^L
\nu\big(\var(f\tc\cF_k)\big) + \frac1L\sum_{k=1}^L \var\big(\nu(f\tc\cF_k)\big)
\,.
\end{equation*}
where $\cF_k$ denote the $\si-$algebra generated by the stick-variables
$\eta_k$, $k=1,\dots,,L$. It is easy to check (see section 3 below) that
the first term can be bounded in terms of $\g(L-1,H)\times \cD(f,f)$.
The main new idea comes in the analysis of the
second term and consists in introducing the stochastic symmetric
operator
\begin{equation*}
Pf =
\frac1L\sum_{k=1}^L \nu(f\tc\cF_k)\,,
\end{equation*}
and observing that for any mean zero function $f$ the following
identity holds true:
\begin{equation*}
\frac1L\sum_{k=1}^L
\var\big(\nu(f\tc\cF_k)\big) = \nu\bigl(fPf\bigr)\,.
\end{equation*}
Thus
$$
\var(f) - \frac1L\sum_{k=1}^L \var\big(\nu(f\tc\cF_k)\big) =
\var(f) - \nu\bigl(fPf\bigr) = \nu\bigl(f(\id-P)f\bigr)\,.
$$ 
so that one is left with the problem of establishing an estimate from
below on the spectral gap of $P$ which is sharp enough to allow a
successful iteration in $L$ for $\g(L,H)$.  The key point now is that,
because of the particular form of $P$ and of the symmetry of the measure
$\nu$, the estimate of the spectral gap of $P$ boils down to the
estimate from below of the spectral gap of a particular one dimensional
random walk that can be described as follows.  Let $n_{\pm}$ be the
minimum and maximum number of particles allowed in a single stick, say
the first one. Then the state space for the random walk is the interval
$[n_-,n_-+1, \dots, n_+]$ and the transition kernel $q(n\to m)$ is given
by $\nu\big(n_1 =m \tc n_2 =n\big)$. It is easy to check that such a
process is ergodic iff $L\ge 3$. The study of its relaxation time
represents in some sense the technical core of the paper and it requires
a rather fine analysis based on equivalence of ensembles type of
results.  \\ 
The rest of the paper is organized as follows.
\begin{enumerate}[(i)]
\item In the next section we define and analyze the one dimensional
random walk mentioned above.
\item In section 3 we prove Theorem \ref{main}.
\item In sections 4,5 we discuss the main applications of Theorem
\ref{main} to quantum Heisenberg XXZ models.
\end{enumerate}

\subsection*{Acknowledgments}
Part of this work was done at the Institute H. Poincar\'e during the
special semester on ``Hydrodynamic limits''. We would like to thank the
organizers F. Golse and S. Olla for their kind invitation and the
stimulating scientific atmosphere there. We are also grateful to
T. Koma, B. Nachtergaele and S. Starr for informing us about their
results prior to publication and to I. Benjamini, N. Berger, C. Hoffman
and E. Mossel for an interesting discussion concerning the mixing time of
asymmetric simple exclusion. Finally we would like to thank F. Cesi for
a very interesting and enlightening discussion which helped us to
clarify a tricky point in our argument.

\newpage

\section{Spectral gap of the long--jump random
walk in a single stick}
In this section we are going to study a one-dimensional process
which plays a key role in the recursive proof of Theorem \ref{main}.

\smallskip

Let $\nu_0$ be the marginal of $\nu$, the canonical measure defined in
(\ref{canon}), on a single stick configuration.  If
$\pi_i:\{0,1\}^\L\to\{0,1\}^{\L_i}$ denotes the canonical projection
onto single stick configurations ($\pi_i \al = \eta_i$) we may write
$\nu_0 = \nu\circ\pi_1^{-1}$.  Let $\bbH$ denote the space
$L^2(\{0,1\}^{\L_1},\nu_0)$ and denote by $\scalar{\cdot}{\cdot}$ the
corresponding scalar product:
$$
\scalar{\varphi}{\psi} = \nu\big((\varphi\circ\pi_1)(\psi\circ\pi_1)\big)
\,,\quad\quad
\varphi,\psi\in\bbH\,.
$$
The $\si$-algebra generated by the single stick variable $\eta_i$ is
denoted by $\cF_i$. The operator $K:\bbH\to\bbH$ is defined by
$$
K\varphi = \nu\big(\varphi\circ\pi_2\tc \cF_1\big)\,,
$$
or equivalently by the bilinear form
$$
\scalar{\varphi}{K\psi} = \nu\big((\varphi\circ\pi_1)(\psi\circ\pi_2)\big)\,.
$$
$K$ is a stochastic, symmetric linear operator on $\bbH$.
The number of particles in a stick is denoted by $n$
and we call $\rho=N/L$ its average value. The centered variable
$n-\r$ is denoted by $\bar n$.
Observe that $\bar n$ is an eigenfunction of $K$ with eigenvalue
$-1/(L-1)$:
\be
K\bar n = -\frac1{L-1}\,\,\bar n\,.
\la{kn}
\end{equation}
Our main result in this section states that
apart from the values $1$ and $-1/(L-1)$ the spectrum of $K$
is concentrated around $0$ within an interval of radius
$O(L^{-1-\d})$ for some $\d>0$, uniformly in $N$, $H$.
\begin{Th}
\la{newteo} There exist constants $\d>0$, $k<\infty$ and $L_0<\infty$
independent of $N,H$ such that if $L\geq L_0$ \be
\big|\scalar{\varphi}{K\varphi}\big|\leq k\,L^{-1-\d}
\scalar{\varphi}{\varphi} \la{theclaim}
\end{equation}
for all $\varphi\in\bbH$ with $\nu_0(\varphi)=0$ and
$\scalar{\varphi}{\bar n}= 0$.
\end{Th}
\proof
Let $\bbH_0$ denote the subspace of $\bbH$ of functions which
only depend on the number of particles. Let also $E$ denote
the orthogonal projection onto $\bbH_0$, i.e.\
if $\cF_0$ denotes the $\si$-algebra generated by the random variable
$n$ we have
$$
E\varphi = \nu_0\big(\varphi\tc \cF_0\big)\,.
$$
A simple computation now shows that $K$ commutes with $E$
and
$$
K\varphi = K E \varphi\,,\quad\quad\varphi\in\bbH\,.
$$
These observations prove that
$\scalar{\varphi}{K\varphi} = \scalar{E\varphi}{KE\varphi}$
and in order to
prove the claim (\ref{theclaim}) we can restrict to $\varphi\in\bbH_0$.

Let us introduce some handy notations. We write $\nu(n)$ for the
probability that in one given stick there are $n$ particles, $\nu(n\cap
m)$ for the probability that in two given sticks there are $n$ and $m$
particles respectively, and $\nu(n\tc m)$ for the probability that in
one given stick there are $n$ particles conditioned to the event that in
a different given stick there are $m$ particles.  With these symbols we
have
$$
K\varphi(n) = \sum_m \nu(m\tc n)\varphi(m)\,\quad\quad \varphi\in\bbH_0\,.
$$
The operator $K$ then describes a random walk on the integers
with transition probabilities
$\text{Prob}\,(n\to m) = \nu(m|n)$.
The desired estimate (\ref{theclaim}) can be written as
\begin{gather}
\sum_{n,m}\nu(m)\nu(n) \varphi(m)Q(m,n)\varphi(n) \,\leq k\,L^{-1-\delta}\,\scalar{\varphi}{\varphi}\,,
\label{teo2} \\
Q(m,n):= \frac{\nu(n\tc m)}{\nu(n)} -1 \,
\label{teo3}
\end{gather}
for all $\varphi\in\bbH_0$ such that $\nu_0(\varphi)=0$ and
$\scalar{\varphi}{\bar n}= 0$.
The rest of this section is concerned with the proof of (\ref{teo2}).

\smallskip

The idea 
is to split the sum in (\ref{teo2}) in a region of typical values of
$m,n$ where things are controlled by a careful expansion and a region of
atypical values whose contribution is shown to be negligible by tail
estimates.  Unfortunately the definition of typical and atypical values
of $n,m$ strongly depends on the value of the particle density $\rho$
and we will be forced to distinguish between two cases, conventionally
denoted {\tt large density} and {\tt small density} case, depending
whether $\rho \ge L^{-\frac{3}{4}}$ or $\rho <
L^{-\frac{3}{4}}$. Before entering into the details of the proof we
first need to establish some preliminary useful bounds.

\subsection{Technical bounds}
The grand canonical
measure $\mu=\mu^{\l(\rho)}$ with density $\rho=N/L$
is the product measure on $\L$ (see (\ref{mula})) with
\be
\mu^{\l(\rho)}(\al)
=\prod_{i=1}^L \prod_{h=1}^H
\frac{q^{2(h-\l(\rho))\al_{(i,h)}}}{1+q^{2(h-\lambda(\rho))}}
\la{mula2}
\end{equation}
Here the parameter $\l(\rho)\in\bbR$, often called the chemical
potential, is such that the average number of particles in any given
stick is equal to $\rho$. The variance $\mu(\bar n^2)$ of the number of
particles in a stick is denoted by $\si^2$. Simple computations - as in
\cite{CapMar}, Lemma 3.2 - show that for every $q\in(0,1)$, there exists
$k<\infty$ such that
$$
k^{-1}(\rho\wedge 1)\leq \si^2 \leq k\,(\rho\wedge 1)\,.
$$
Similarly $\mu(\bar n^4)\leq k(\rho\wedge 1)$.
As a rule, here and throughout the rest
of this section the letter $k$ will be used to denote
a finite constant whose value may change from line to line.
What is essential is that it only depends on $q$ and is uniform in
all other parameters: $L$,$H$, and $N$.

\smallskip

We introduce the characteristic function
$$
F(t)=\mu\big(\nep{i\frac{t}{\si\sqrt{L}}\bar n}\big)\,.
$$
We shall rely on the following simple estimate.
\begin{Le}
\la{gauss}
For all $t\in [-\pi\si\sqrt{L},\pi\si\sqrt{L}]$ 
\be
|F(t)|\leq \nep{-k \frac{t^2}{L}}
\la{gauss1}
\end{equation}
\end{Le}
\proof
Writing
$$
n = \sum_h \al_h \,,
$$
where we drop the horizontal stick label, we have
$$
|F(t)| = \prod_h g_h(t)\,,\quad\quad
g_h(t):=\big|\mu\big(\nep{i\frac{t}{\si\sqrt{L}}\al_h}\big)\big|\,.
$$
Using the inequality $x\leq\nep{(x^2-1)/2}$, $x\in(0,1)$, we have
\begin{align*}
g_h(t)&\leq \exp{\big(\frac12(g_h(t)^2-1)\big)} \\
& = \exp{\big\{-\frac12\var_\mu\big[\cos{(t\al_h/{\si\sqrt{L}})}\big]
-\frac12\var_\mu\big[\sin{(t\al_h/{\si\sqrt{L}})}\big]\,\big\}}
\end{align*}
Define $\si^2_h = \var_\mu(\al_h)$ and compute
\begin{gather*}
\var_\mu\big[\cos{(t\al_h/{\si\sqrt{L}})}\big] = \si^2_h
\big(\cos{(t/{\si\sqrt{L}})} - 1\big)^2\\
\var_\mu\big[\sin{(t\al_h/{\si\sqrt{L}})}\big] = \si^2_h
\big(\sin{(t/{\si\sqrt{L}})}\big)^2\,.
\end{gather*}
It follows
$$
g_h(t)\leq\exp{\big\{- \si^2_h
(1 - \cos{(t/{\si\sqrt{L}})})\big\}}\leq
\exp{\big(-kt^2\si^2_h/\si^2L\big)}\,,
$$
where we use the inequality $1-\cos s \geq k s^2$, for some $k>0$ and all
$s\in[-\pi,\pi]$. The lemma now follows from $\si^2=\sum_h \si^2_h$.
\qed

\\
It is useful to compare the
canonical measure $\nu$, given again by (\ref{canon}), with the
grand canonical probability $\mu$. In particular, it follows from
Lemma \ref{gauss} - as shown in \cite{CapMar}, Proposition 3.8 -
that there exists $k<\infty$ independent of the density $\r$ such that
\be
\nu(f)\leq k\,\mu(f)
\la{comparison}
\end{equation}
for any function $f\ge 0$ depending at most on $\frac{L}{2}$ variables
$\{\h_i\}$.  We shall need the following estimates on the tails of our
distributions.
\begin{Le}
\la{tail} 
There exist constants $a>0$ and $k<\infty$ depending only on
$q$ such that
\begin{align}
\nu(n\tc m)\nu(m\tc n) & \le k\, \nep{ -\,a\,[(n-\rho)^2 + (m-\rho)^2]}
\la{tail1}\\
\nu(n\tc m)\nu(m\tc n) &\le k\, \rho^{n+m}
\la{2}
\end{align}
\end{Le}
\proof
We start with the proof of (\ref{tail1}).
Given $n\in\bbN$ we denote by $\l(n)$ the chemical
potential such $\mu^{\l(n)}(n_1)=n$. We write simply $\l$
for the chemical potential $\l(\rho)$.
We have
\begin{equation}
\nu(n)  \le k\,\mu^\l(n) = k\,
\nep{-V_\l(n)} \mu^{\l(n)}(n) \le k\, \nep{-V_\l(n)}
\la{su}
\end{equation}
where
$$
V_\l(n) = c\,(\l(n)-\l)n -\log\(\frac{Z^{\l(n)}}{Z^\l}\)   \,,
\quad c:=-\log q\,,
$$
and $Z^\l$ denotes the partition function
$$
\prod_{h=1}^H {\(1+q^{2(h-\lambda)}\)}\,.
$$
Now $V_\l(\rho)=0$ and we write
$$
V_\l(n) = \int_{\rho}^n ds \, \frac{d}{ds} V_\l(s)\,.
$$
A standard
computation gives that $\frac{d}{ds} V_\l(s) =
c\,(\l(s)-\l)$ and
$$
\frac{d^2}{ds^2} V_\l(s) = c\,\frac{d}{ds}\l(s)=
\frac1{\mu^{\l(s)}(\bar n ^2)} \geq a
$$
for some $a>0$ independent of $s$.  Therefore
$$
V_\l(n) \geq a\, (n-\rho)^2\,.
$$
By (\ref{su}) we have the tail estimate
\be
\nu(n) \leq k \nep{-a\,(n-\rho)^2}\,.
\la{tail0}
\end{equation}
We now apply the above bound to each term $\nu(n\tc m)$. In this case
the starting chemical potential $\l$ corresponds to a density
$\rho' = \frac{N-m}{L-1}$ so that (\ref{tail0}) yields
$$
\nu(n\tc m) \le k\, \nep{-a\,\(n-\frac{N-m}{L-1}\)^2}
$$
Our original claim (\ref{tail1}) now follows from the uniform (in $N,L$)
convexity of the function
$$
  G(x,y) := \(x-\frac{N-y}{L-1}\)^2 + \(y-\frac{N-x}{L-1}\)^2
$$
around the unique minimum $x= \rho$, $y=\rho$.

\\
To prove (\ref{2}) we use again the
bound (\ref{comparison}) and estimate
$$
\nu(n\tc m) \le k \mu^{\l(\frac{N-m}{L-1})}(n)\le
k \(\frac{N-m}{L-1}\)^n \le k \rho^n\,.
$$
\qed

\\
We are in a position to complete the proof of Theorem \ref{newteo}.  As
anticipated we will need to distinguish between ``large'' and ``small''
values of the density $\rho$.

\subsection{The large density case}
Here we assume $\rho\geq L^{-\frac{3}{4}}$.
Define the set
\be
\cB = \{(n,m):\;|\bar n| + |\bar m| \leq B\log L\}\,,
\la{b}
\end{equation}
with $B<\infty$ a constant to be fixed later.
We rewrite the sum in (\ref{teo2}) as
\be
\sum_{(n,m)\in\cB}\nu(n)\nu(m)
\varphi(n)Q(m,n)\varphi(m) + \sum_{(n,m)\notin\cB}\nu(n)\nu(m)
\varphi(n)Q(m,n)\varphi(m)\,.
\la{teo5}
\end{equation}
The second term in (\ref{teo5}) is estimated with the help of
Lemma \ref{tail}. To see this use Schwarz' inequality
to write
\begin{align*}
&\sum_{(n,m)\notin\cB}
\nu(m)\nu(n)\varphi(m) Q(m,n) \varphi(n)  \\ &\le
\scalar{\varphi}{\varphi}\Bigl[\sum_{(n,m)\notin\cB}
\nu(m)\nu(n)Q(m,n)^2 \Bigr]^{\frac12} \\
&\le  \scalar{\varphi}{\varphi}\Bigl[\sum_{(n,m)\notin\cB}
\nu(n)\nu(m) 2\,[1 + \frac{\nu(n\tc m)^2}{\nu(n)^2}] \, \Bigr]^{\frac12}\,.
\end{align*}
Since
$$
\nu(n)\nu(m)\frac{\nu(n\tc m)^2}{\nu(n)^2}
= \nu(n\tc m) \nu(m \tc n)\,,
$$
an application of
(\ref{tail1}) and (\ref{tail0}) gives
\be
\sum_{n,m\notin\cB}\nu(n)\nu(m)
\varphi(n)Q(m,n)\varphi(m) \leq \,k\, L^{-2}\,\scalar{\varphi}{\varphi}
\,,
\la{sor}
\end{equation}
provided $B$ is sufficiently large.

\smallskip

The key step in the proof of Theorem \ref{newteo} in the case of
large density will be the following expansion.
\begin{Le}
\la{expa}
For all $q\in(0,1)$, $B<\infty$,
there exist constants $k<\infty$, $\zeta > 0$ such that
if $\rho\geq L^{-\frac34}$
then
\be
Q(m,n) = -\frac{\bar n \bar m}{\si^2 L}  +
\cR(n,m)
\la{expa1}
\end{equation}
with the remainder $\cR$ satisfying
\be
\Bigl[\,
\sumtwo{n,m}{|\bar n|+|\bar m|\leq B\log{L}} \nu(n)\nu(m) |\cR(n,m)|^2
\, \Bigr]^{1/2} \le k L^{-1-\zeta}\,.
\la{expa2}
\end{equation}
\end{Le}
\proof
We write
\begin{align}
\nu(n) &=\,
        \frac{\mu(n)}{2\pi\si\sqrt{L}\mu(N)}\int dt \, F(t)^{L-1}
         \nep{i\frac{t}{\si\sqrt{L}}\bar n} \nonumber \\
\nu(n\cap m) &=\,
        \frac{\mu(n)\mu(m)}{2\pi\si\sqrt{L}\mu(N)}\int dt \, F(t)^{L-2}
         \nep{i\frac{t}{\si\sqrt{L}}[\bar n +\bar m]} \nonumber
\end{align}
where all the integrals are over the interval
$[-\pi\si\sqrt{L},\pi\si\sqrt{L}]$. Since
$$
2\pi\si\sqrt{L}\mu(N) = \int dt \, F(t)^{L}\,
$$
we can write
$$
Q(m,n) = \frac{\nu(n\cap m)-\nu(n)\nu(m)}{\nu(n)\nu(m)} =
\frac{\rm NUM}{\rm DEN}
$$
with
\begin{gather*}
\text{NUM} :=
\int dt \, F(t)^{L-2}
         \nep{i\frac{t}{\si\sqrt{L}}[\bar n +\bar m]}
\int dt' \, F(t')^{L} \\
- \int dt \, F(t)^{L-1}  \nep{i\frac{t}{\si\sqrt{L}}\bar n }
         \int dt' \, F(t')^{L-1} \nep{i\frac{t'}{\si\sqrt{L}}\bar m}
\end{gather*}
and
$$
\text{DEN} :=
\int dt \, F(t)^{L-1} \nep{i\frac{t}{\si\sqrt{L}}\bar n }
         \int dt' \, F(t')^{L-1}\nep{i\frac{t'}{\si\sqrt{L}}\bar m}
$$
Notice that, because of the Gaussian upper bound of Lemma \ref{gauss}
we have $|F(t)|^L \leq \nep{- at^2}$ and only the region $|t| \le
k\log L$ (for some large but fixed $k$) will have to be taken care of.
%
In order to be precise about the nature of the various error terms, in
what follows we will denote by $\e(L)$ a generic term which, upon
multiplication by $L$, still goes to zero (as $L\to \infty$) as an inverse
power of $L$ uniformly in the range of $|t|\leq k\log{L}$.
We first observe that $F(t) = 1-\frac{t^2}{2L} + \e(L)$.
Indeed by expanding $F$ around $t=0$ the third order error term
is bounded from above by
$$
k\frac{|t|^3}{(\si^2L)^{3/2}}\mu( |\bar n|^3)
\le k \frac{|t|^3}{\si L^{3/2}} \le \e(L)
$$
where we use the bound
$$
\mu( |\bar n|^3)\leq k\,\si \,\mu(\bar n^4)^\frac12\leq k\,\si^2
$$
together with
$\si^2\geq k(\rho\wedge 1)$ and $\rho\geq L^{-\frac34}$.
This implies $F(t)^{-1}=1+\frac{t^2}{2L} + \e(L)$ and
$F(t)^{-2}=1+\frac{t^2}{L} + \e(L)$. Then
if we write $\nep{i\frac{t}{\si\sqrt{L}}\bar n} = 1+\d_n(t)$ we have
$$
\text{NUM} = \int dt  F(t)^{L} \d_n(t)\d_m(t)
            \int dt'  F(t')^{L}   -
           \int dt  F(t)^{L}\d_n(t) \int dt' F(t')^{L}\d_m(t') +\e(L)\,.
$$
If we modify further and define
$$
\hat \d_n = \d_n -
i\frac{t}{\si\sqrt{L}}\bar n\,,\quad\quad I(L)=\int dt  F(t)^{L}\,t
$$
we then have,
\begin{gather*}
\text{NUM} = - \frac{\bar n \bar m}{\si^2 L}\Big(\int dt \, t^2 \,F(t)^{L}
\int dt' \, F(t')^{L} - I(L)^2\Big)\\
+ \int dt \, F(t)^{L}\hat\d_n(t)\,\hat\d_m(t)
            \int dt' \, F(t')^{L}   -
   \int dt \, F(t)^{L}\hat\d_n(t) \int dt' \, F(t')^{L}\hat\d_m(t') \\
  +\, i\,\frac{\bar m}{\si\sqrt{L}}
  \,\int dt \, F(t)^{L}\,t\,\hat\d_n(t)\int dt' \, F(t')^{L}
  + \,i\,\frac{\bar n}{\si\sqrt{L}}
  \,\int dt \, F(t)^{L}\,t\,\hat\d_m(t)\int dt' \, F(t')^{L} \\
  - \,i\,\frac{\bar m}{\si\sqrt{L}}
      \int dt \, F(t)^{L}\hat\d_n(t)
       \int dt' \, F(t')^{L}\,t'
   - \,i\,\frac{\bar n}{\si\sqrt{L}}
     \int dt \, F(t)^{L}\hat\d_m(t)
       \int dt' \, F(t')^{L} \,t'
                          \,+ \,\e(L)\,.
\end{gather*}
Now observe that $\frac1LI(L)=\e(L)$ and
$$
\frac{1}{L}\,\int dt \, t^2 \,F(t)^{L}
\int dt' \, F(t')^{L}
= \frac{{2\pi}}{L} + \e(L)\,.
$$
Using also
$|\hat \d_n(t)| \le k \frac{\bar n^2 t^2}{\si^2 L} $ it follows that
$$
\text{NUM} = - 2\pi\,\frac{\bar n \bar m}{\si^2 L}
+ R(n,m) +\e(L)
$$
with
\begin{equation*}
|R(n,m)| \le k \Bigl[\,\frac{|\bar n|\,|\bar m| }{\si^2}\,\e(L) +
 \frac{\bar n^2\bar m^2 }{(\si^2 L)^2} +
\frac{|\bar n|\,\bar m^2 + |\bar m|\,\bar n^2}{(\si^2 L)^{3/2}}
\,\Bigr]
\la{expa45}
\end{equation*}
What is crucial for us is that
\be
\Bigl[\,
\sumtwo{n,m}{|\bar n|+\bar m|\leq B\log{L}}
\nu(n)\nu(m)|R(n,m)|^2 \, \Bigr]^{1/2} \le \e(L)
\la{expa5}
\end{equation}
provided that $\rho \ge L^{-\frac34}$.
The above estimate actually holds without
the restriction $|\bar n|+|\bar m|\leq B\log{L}$ as it is easily
seen using $\nu(\bar n^4)\leq k\mu(\bar n^4) \leq k\si^2$ and
the bound $(\rho\wedge 1)\leq k \si^2$.
On the other hand if we repeat the reasoning for
the denominator $\text{DEN}$ we obtain
\begin{equation*}
\text{DEN} = 2\pi + \hat R(n,m)
\la{expa46}
\end{equation*}
with a remainder $\hat R(n,m)$ satisfying
\be
\suptwo{n,m:}{|\bar n|+|\bar m|\leq B\log{L}} |\hat R(n,m)| \leq
k L^{-\zeta}
\la{expa47}
\end{equation}
for some $\zeta > 0$ whenever $\rho\geq L^{-\frac34}$.
In conclusion $Q(m,n)$ has been written as in (\ref{expa1}) and
(\ref{expa5})--(\ref{expa47}) imply (\ref{expa2}). \qed

\bigskip

We are now able to finish the proof of Theorem \ref{newteo} in the case of
large densities. %
Observe that
$$
\sum_{(n,m)\in\cB}\nu(n)\nu(m)
\varphi(n)\varphi(m)\bar n \bar m =
\scalar{\varphi}{\bar n}^2 - \sum_{(n,m)\notin\cB}\nu(n)\nu(m)
\varphi(n)\varphi(m)\bar n \bar m \,.
$$
By assumption $\scalar{\varphi}{\bar n}=0$. Then
\begin{align*}
\Big|
\sum_{(n,m)\in\cB}&\nu(n)\nu(m)
\varphi(n)\varphi(m)\bar n \bar m
\,\Big| \\
&\leq \Big[
\sum_{(n,m)\notin\cB}\nu(n)\nu(m)
\bar n^2 \bar m^2
\Big]^\frac12
\,\scalar{\varphi}{\varphi}
\leq k L^{-2}\,\scalar{\varphi}{\varphi}
\end{align*}
where the last estimate can be easily obtained
from (\ref{tail0}).
By Lemma \ref{expa} we then have
\begin{align*}
\Big|
\sum_{(n,m)\in\cB}&\nu(n)\nu(m)
\varphi(n)Q(n,m)\varphi(m)
\,\Big| \\
&\leq \scalar{\varphi}{\varphi}
\Big\{
\,k\,L^{-2} \,+\, \Big[\,
\sum_{(n,m)\in\cB} \nu(n)\nu(m) |\cR(n,m)|^2
\, \Big]^{1/2}
\Big\}
\leq \,k\, L^{-1-\zeta}\,\scalar{\varphi}{\varphi}
\end{align*}
Together with (\ref{sor}) these bounds imply Theorem \ref{newteo}
when $\rho\geq L^{-\frac34}$.

\subsection{The case of small density}
When $\rho \le L^{-\frac34}$ the strategy for
the proof of (\ref{teo2}) has to be slightly modified since various
previous technical estimates are no longer valid. We can however take advantage
of the thinner tails of the distribution of the number of particles. In
this respect we observe that
\be
\Bigl[\sumtwo{n,m}{n+m\geq 2\,,\:nm
\neq 1} \nu(n)\nu(m) \,Q(n,m)^2 \, \Bigr]^{\frac12} \le k L^{-1-\zeta}
\la{small}
\end{equation}
with $\zeta = 1/8$.
Indeed, thanks to (\ref{2})
$$
\Bigl[\sumtwo{n,m}{n+m\ge 3}
\nu(n)\nu(m) \,Q(n,m)^2 \, \Bigr]^{\frac12} \le k \,\rho^{\frac32}\leq
k\,L^{-1-\zeta}\,.
$$
On the other hand
we can examine explicitly  the case
$n=2, m=0$. We have
$$
Q(2,0) = \frac{1}{\nu(0)} - 1 -
     \frac{\nu(\{n=2\}\cap\{ m\ge 1\})}{\nu(0)\nu(2)} \,.
$$
Using the bounds $\nu(0)\geq 1-k\rho$, $\nu(n)\leq k\rho^n$
and $\nu(\{n=2\}\cap\{ m\ge 1\})\leq k\rho^3$ one has
$$
\nu(2)\nu(0)Q(2,0)^2
\leq k\rho^4\,.
$$
We have proved (\ref{small}).

\smallskip

At this point we may proceed as in (\ref{teo5}) with the choice
$$\cB := \{ n,m: \; n \le 1, \, m\le
1\}\,.$$
The second term in (\ref{teo5}) is
controlled by (\ref{small}).
The first term is given by
$$
\sum_{n\le 1}\sum_{m\le 1} \varphi(n)\varphi(m)\bigl[\, \nu(n\cap
m)-\nu(n)\nu(m)\,\bigr]\,.
$$
We are going to study the eigenvalues of the symmetric $2\times 2$
matrix
$$
\cM(n,m):= \bigl[\, \nu(n\cap
m)-\nu(n)\nu(m)\,\bigr]\,,\quad n,m\in\{0,1\}\,.
$$
For $n=0,1\dots N$, let $M_n$ be the number of sticks
with exactly $n$ particles.
The identities $L=\sum_n M_n$ and $N=\sum_n nM_n$ imply
$$
M_0 = L - N +\sum_{n\ge 2}(n-1)M_n\,,\quad\quad
M_1 = N -\sum_{n\ge 2}nM_n \,.
$$
Denoting by $n_i$ the number of particles
in the stick $\L_i$, we have the estimates
\begin{align*}
\var_\nu(M_0) &\le k \var_\mu\(\sum_i
(n_i-1)\id_{\{n_i\geq 2\}}\) \le kL\rho^2 \\
\var_\nu(M_1) &\le k \var_\mu\(\sum_i n_i\id_{\{n_i\geq
2\}}\) \le kL\rho^2 \,.
\end{align*}
We compute
\begin{align*}
\nu(\{n_1=0\}&\cap\{n_2=0\}) = \frac1{L-1}\sum_{i=2}^L
\nu(\{n_1=0\}\cap\{n_i=0\}) \\
& = \frac1{L-1}\,\nu(\id_{\{n_1=0\}}(M_0-1))
= \frac1{L(L-1)}\,\nu(M_0(M_0-1)) \\
& = \frac1{L(L-1)}\,\var(M_0) +
\frac1{L(L-1)}\,\nu(M_0)^2 - \frac1{L(L-1)}\,\nu(M_0)\,.
\end{align*}
It follows that $\cM(0,0) = \nu(\{n_1=0\}\cap\{n_2=0\}) - \nu(1)\nu(0)$
can be written
$$
\cM(0,0) = - \frac{\nu(0)(1-\nu(0))}{L-1} + \frac1{L(L-1)}\var(M_0)\,.
$$
Similarly
\begin{align*}
\cM(1,1) &= - \frac{\nu(1)(1-\nu(1))}{L-1} + \frac1{L(L-1)}\var(M_1) \\
\cM(0,1) &= \frac{\nu(0)\nu(1)}{L-1} + \frac1{L(L-1)}\,{\rm Cov}_\nu(M_0,M_1)\,.
\end{align*}
Introducing
the matrix
$$
A =
\begin{pmatrix}
-1 & 1\\
1 & -1
\end{pmatrix}\,,
$$
the above computations show that
$$
\cM = \frac{\r}{L-1}\,A \,+\, \widetilde\cM
$$
with a symmetric matrix $\widetilde\cM$ such that
$|\widetilde\cM(n,m)|\leq O(\r^2/L)$, $n,m\in\{0,1\}$.
Then we have
\begin{equation}
\Big|
\sumtwo{n\leq 1}{m\leq 1}\varphi(n)\varphi(m)\cM(n,m)
\Big|
\leq \frac{\r}{L-1}\, [\varphi(0)-\varphi(1)]^2 + k\frac{\r^2}{L}
\,[\varphi(1)^2 + \varphi(0)^2]
\la{ba}
\end{equation}
Observe that by the orthogonality $\scalar{\varphi}{\bar n}=0$
and Schwarz inequality
\begin{align*}
\Big|\sumtwo{n\leq 1}{m\leq 1}&\nu(n)\nu(m)
\varphi(n)\varphi(m)\bar n \bar m\Big|\\
&\leq \scalar{\varphi}{\varphi}
\,\Big[\sumtwo{n+m\geq 2}{nm\neq 1}\nu(n)\nu(m)
\bar n^2 \bar m^2\Big]^\frac12
\leq k\,\r^\frac32\, \scalar{\varphi}{\varphi}
\end{align*}
On the other hand
$$
\sumtwo{n\leq 1}{m\leq 1}\nu(n)\nu(m)
\varphi(n)\varphi(m)\bar n
\bar m = \r^2[\varphi(0)-\varphi(1)]^2
+ O(\r^3)[\varphi(1)^2 + \varphi(0)^2]\,.
$$
Collecting these estimates we arrive at
$$
[\varphi(0)-\varphi(1)]^2\leq
k \,\r^{-\frac12} \,\scalar{\varphi}{\varphi} +
k\,\r\,[\varphi(1)^2 + \varphi(0)^2]
$$
Using also the bound
$$
[\varphi(1)^2 + \varphi(0)^2]\leq k\r^{-1}\scalar{\varphi}{\varphi}
$$
and going back to (\ref{ba}) we finally obtain
$$
\Big|
\sumtwo{n\leq 1}{m\leq 1}\varphi(n)\varphi(m)\cM(n,m)
\Big|\leq k \frac{\r^\frac12}{L}  \,\scalar{\varphi}{\varphi}
\leq k L^{-1-\zeta}  \,\scalar{\varphi}{\varphi}\,.
$$
This ends the proof of Theorem \ref{newteo}

\newpage
\section{Proof of Theorem \ref{main}}
%
Recall the definition of the stochastic operator $K$ introduced at the
beginning of the previous section. We set \be w(L,H) = \sup_{N}
\big[\gap(\id - K)\big]^{-1} \la{wlh}
\end{equation}
where
\be
\gap(\id-K) = \inftwo{\varphi\in\bbH:}{\nu_0(\varphi)=0}\,
\frac{\scalar{\varphi}{(\id-K)\varphi}}{\scalar{\varphi}{\varphi}} \,.
\la{gapk}
\end{equation}
Define also
\be
\tg(L,H) = \sup_{N} \sup_{f\in L^2(\nu)}\,
\frac{\var(f)}{\tilde\cD(f,f)}
\la{tgammalh}
\end{equation}
where we introduce the modified Dirichlet form
\be
\tilde\cD(f,f) = \frac1L \sum_{i} \sum_{j\neq i} D_{ij}(f)\,.
\la{tcd}
\end{equation}
Clearly $\cD(f,f)\geq \tilde\cD(f,f)$ for all $f$ and $\g(L,H)\leq
\tg(L,H)$. Note also that $\tilde\cD(f,f)$ is still ergodic:
$\tilde\cD(f,f)=0$ implies $f={\rm const.}$ for all $L\geq 3$. More precisely,
using e.g.\ the method of Lemma 2.6 in \cite{CapMar} one can easily
prove that there exists $k=k(q)<\infty$ such that for every $L\geq 3$
one has $\cD(f,f)\leq k\tilde\cD(f,f)$ for all $f\in L^2(\nu)$. The key
step in the proof of Theorem \ref{main} is represented by the following
proposition whose proof is postponed to the end of the section.
\begin{Pro}
\la{pro_iter}
There exists $L_0<\infty$ such that
for any $L\geq L_0$ and any $H\geq 1$
\be
\tg(L,H) \leq \big[1\vee w(L,H)\big]\, \tg(L-1,H)
\la{iter}
\end{equation}
\end{Pro}

\\
To complete the proof
of Theorem \ref{main} we need an estimate on $w(L,H)$.
In view of Theorem \ref{newteo} we know that for
every $\varphi\in\bbH$, such that $\nu_0(\varphi)=0$ we
have
\be
\scalar{\varphi}{(\id-K)\varphi} \geq \big(1-k\,L^{-1-\d}\big)
\scalar{\varphi}{\varphi}
\la{oclaim}
\end{equation}
for $L\geq L_0$
with uniform constants $\d>0$ and $L_0,k<\infty$.
It is then immediate to deduce $w(L,H)\leq 1 + k L^{-1-\d}$ and therefore
\be
\prod_{L=L_0}^\infty w(L,H) \leq C\,,
\la{claim}
\end{equation}
with some uniform constant $C<\infty$.
Then Theorem \ref{main} follows from Proposition \ref{pro_iter} and the bound
$$
\sup_H\tg(L_0,H)<\infty\,,\quad\quad L_0 \geq 3\,.
$$
The latter 
is easily deduced e.g.\ from results in \cite{CapMar}.

\subsection{Proof of Proposition \ref{pro_iter}}
Let $\cF_k$ denote the $\si-$algebra generated by the stick-variables
$\eta_k$, $k=1,\dots,,L$. For any $f$ one has the decomposition
\be
\var(f) = \frac1L\sum_{k=1}^L
\Big\{
\nu\big(\var(f\tc\cF_k)\big) + \var\big(\nu(f\tc\cF_k)\big)
\Big\}\,.
\la{iter10}
\end{equation}
We first establish the estimate
\be
\frac1L\sum_{k=1}^L
\nu\big(\var(f\tc\cF_k)\big)
\leq \frac{L-2}{L-1} \,\tg(L-1,H) \tilde\cD(f,f)\,.
\la{iter20}
\end{equation}
By definition of $\tg$ for every $k$ one has
$$
\nu\big(\var(f\tc\cF_k)\big) \leq \tg(L-1,H)
\,\frac1{L-1}\sum_{i\neq k}\sum_{j\neq i,k} D_{ij}(f)\,.
$$
Summing over $k$ and using
$$
\sum_k \sum_{i\neq k}\sum_{j\neq i,k} D_{ij}(f) =
(L-2) \sum_{i}\sum_{j\neq i} D_{ij}(f) = L(L-2) \tilde\cD(f,f)
$$
we obtain (\ref{iter20}).

\smallskip

We turn to an estimate on the second term in (\ref{iter10}).
Consider the non-negative stochastic operator $P:L^2(\nu)\to
L^2(\nu)$ defined by
\be
Pf =
\frac1L\sum_{k=1}^L \nu(f\tc\cF_k)\,.
\la{iter30}
\end{equation}
We may assume without loss that $\nu(f)=0$. Then we have the
identity
\be
\frac1L\sum_{k=1}^L
\var\big(\nu(f\tc\cF_k)\big) = \nu(fPf)\,.
\la{iter50}
\end{equation}
From (\ref{iter10}), (\ref{iter20})
and (\ref{iter50}) we obtain
\be
\nu(f(\id - P)f) \leq \frac{L-2}{L-1} \,\tg(L-1,H) \tilde\cD(f,f)\,.
\la{iter60}
\end{equation}
In order to estimate from below the left hand side in (\ref{iter60})
we are going to prove the bound
\be
\gap(\id-P) \geq \frac{L-2}{L-1}\,\big[1\wedge \gap(\id-K)\big]\,.
\la{claimo}
\end{equation}
Here $\gap(\id-P)$ stands for the smallest nonzero eigenvalue of $\id - P$.
Note that (\ref{claimo}) and (\ref{iter60}) immediately imply the
proposition.
\smallskip

Take $f\in L^2(\nu)$ such that $\nu(f)=0$ and
$$
Pf=\l f\,,\quad\quad \l> 0\,.
$$
Then $f$ is of the form
\be
f(\eta) = \sum_{\ell}\varphi_\ell(\eta_\ell)\,,
\la{eq1}
\end{equation}
with $\varphi_\ell\in\bbH$ and we have the identity
\be
\l \sum_{\ell}\varphi_\ell(\eta_\ell)
= \frac1L \sum_{\ell} \sum_{k} \nu(\varphi_k(\eta_k)\tc\eta_\ell)\,.
\la{id1}
\end{equation}
Define now the function $\Phi_f\in\bbH$:
$$
\Phi_f = \sum_\ell \varphi_\ell\,.
$$
Taking conditional expectation with respect to $\cF_j$
in (\ref{id1}) a simple computation yields
$$
\l K \Phi_f + \l(\id-K)\varphi_j
= \frac{L-1}{L} K^2\Phi_f + \frac1L K(2-K)\Phi_f
+ \frac1L(\id-K)^2\varphi_j\,.
$$
Summing over $j$ and factorizing we obtain
\be
\left[K - \frac{\l L - 1}{L-1}\right] \left[K + \frac{1}{L-1}\right]
\Phi_f = 0
\la{id}
\end{equation}
The above identity says that if $(K+\frac{1}{L-1})\Phi_f\neq 0$ then
$\mu := \frac{\l L - 1}{L-1}$ is in the spectrum of $K$. In this case
then
$$
1-\l = \frac{L-1}L (1-\mu)\geq \frac{L-1}{L} \gap(\id-K) \geq
\frac{L-2}{L-1} \gap(\id-K)\,.
$$
We have to study the eigenfunctions $f$ of $P$ such that
the corresponding $\Phi_f$ satisfies
\be
\left[K + \frac{1}{L-1}\right]
\Phi_f = 0 \,.
\la{ide}
\end{equation}
Then Theorem \ref{newteo} implies that
$\Phi_f = A\bar n$ for some $A\in\bbR$.
On the other hand
every $\varphi_\ell$ can be decomposed
as
$$
\varphi_\ell = a_\ell \bar n + \hat\varphi_\ell
$$
with $a_\ell=\scalar{\varphi_\ell}{\bar n}$ and
$\scalar{\hat\varphi_\ell}{\bar n}=0$.
Note that in view of (\ref{eq1}) and
the conservation law $\sum_\ell \bar n(\eta_\ell) = 0$, 
there is no restriction
in assuming $\sum_\ell a_\ell = 0$.
Therefore
$$A \scalar{\bar n}{\bar n}
= \scalar{\Phi_f}{\bar n} = \scalar{\bar n}{\bar n}
\sum_\ell a_\ell = 0
$$
Thus if $\Phi_f$ solves (\ref{ide}) then necessarily $\Phi_f=0$.
It remains to study the latter case in detail.

\smallskip
Typical examples for which $\Phi_f=0$ are obtained by
choosing each $\varphi_\ell$ proportional to $\bar n$.
Call $\cA$ the class of all $f\in L^2(\nu)$
of the form (\ref{eq1}) with $\varphi_\ell= a_\ell \,\bar n$
with arbitrary $\underline a :=\{a_1,\dots,a_L\}\in \bbC^L$. A simple
computation gives
\be
(\id-P)f = \frac{L-2}{L-1} f\,,\quad \quad f\in\cA
\la{eq2}
\end{equation}
By (\ref{eq2}) we may
then restrict to the orthogonal complement $\cA^{\perp}$ to
prove our claim (\ref{claimo}).
We are going to prove that there exists a finite $L_0=L_0(q)$ such that
if $f\in \cA^{\perp}$ and $\Phi_f=0$
\be
(1-\l)\geq \frac{L-2}{L-1}\,,\quad\quad \forall L\ge L_0.
\la{eq3}
\end{equation}
Let us first check that if $f\in \cA^{\perp}$ then the corresponding
$\varphi_\ell$ are all orthogonal (in $\bbH$) to the number of particles:
\be
\scalar{\varphi_\ell}{\bar n}
= 0\,,\quad\quad \ell =1,\dots,L
\la{eq4}
\end{equation}
Indeed $f\in \cA^{\perp}$ means
\begin{align*}
0 &= \sum_{k,\ell}a_k \nu(\varphi_\ell(\eta_\ell) \bar n_k) =
\sum_k a_k \Big[\sum_{\ell\neq k} \scalar{\varphi_\ell}{K \bar n}
+ \scalar{\varphi_k}{\bar n}\Big] \\
& = \frac{1}{L-1} \sum_k a_k \Big[ - \scalar{\Phi_f}{\bar n}
  + L\,\scalar{\varphi_k}{\bar n}\Big]  =
\frac{L}{L-1} \sum_k a_k \,\scalar{\varphi_k}{\bar n}
\end{align*}
Here we are using (\ref{kn}) and $\Phi_f = 0$.  The above identity
implies (\ref{eq4}) in view of the arbitrariness of $\underline a$.
Using again $\Phi_f=0$ we compute
\begin{align*}
\nu(f^2) &=
\sum_{k,\ell} \nu(\varphi_k(\eta_k)\varphi_\ell(\eta_\ell))
= \sum_{k}\Big[ \sum_{\ell\neq k} \scalar {\varphi_k}{K\varphi_\ell}
+ \scalar{\varphi_k}{\varphi_k}\Big]
\\
& =  \sum_k \Big[\scalar {\varphi_k}{K\Phi_f} +
 \scalar{\varphi_k}{(\id-K)\varphi_k}\Big] =
\sum_k\scalar{\varphi_k}{(\id-K)\varphi_k}
\end{align*}
A similar computation yields
$$
\nu(f(\id-P)f) = \frac1L \sum_{k}\scalar{\varphi_k}{(\id-K)
((L-1)\id+K)\varphi_k}
$$
Writing $\tilde\varphi_\ell = (\id-K)^{1/2}\varphi_\ell$
one has
$$
(1-\l)\sum_\ell\scalar{\tilde\varphi_\ell}{\tilde\varphi_\ell}
= \frac1L\sum_\ell\scalar{\tilde\varphi_\ell}
{((L-1)\id +K)\tilde\varphi_\ell}\,.
$$
Now observe that $\scalar{\tilde\varphi_\ell}{\bar n}=0$. This
follows from (\ref{eq4}), (\ref{kn}) and the self adjointness of $K$.
By Theorem \ref{newteo} we then infer
$$
\big|
\scalar{\tilde\varphi_\ell}{ K \tilde\varphi_\ell}\big|
\leq k \,L^{-1-\d} \,\scalar{\tilde\varphi_\ell}{\tilde\varphi_\ell}
$$
for some uniform constants $\d>0$, $k<\infty$. We conclude that
$$
(1-\l) \geq \frac{L-1}{L} - \frac{k}{L^{2+\d}} \geq \frac{L-2}{L-1}
$$
for $L$ large. This finishes the proof of (\ref{claimo}).
\qed

\bigskip
\subsection{A remark on Bernoulli--Laplace model of diffusion}
We observe that the strategy of the above proof
may be used to compute in a simple way
the spectral gap for the so--called Bernoulli--Laplace process, see
\cite{DSh}, \cite{DS} and references therein.
The latter can be seen as an exclusion process on a complete graph:
there are
$L$ sites with exchanges allowed between any couple of
sites $(i,j)$, $i,j=1,\dots,L$ and
with uniform rates. The configuration space is
$\Om_0:=\{0,1\}^L$ and the measure $\nu$ is the product Bernoulli
measure on $\Om_0$ conditioned to the event
$\sum_i \a_i = N$. The Dirichlet form is then defined by (\ref{tcd})
with $D_{ij}$ replaced by
\begin{gather}
E_{ij}(f) = \frac12
\nu\big[(\grad_{ij}f)^2\big]
\la{eij}\\
\grad_{ij}f(\a) = f(\a^{ij})-f(\a)\,,\quad
(\a^{ij})_k = \begin{cases}
    \a_k & k\neq i,j\\
    \a_i & k=j \\
    \a_j & k=i
             \end{cases}
\end{gather}
We are going to show that $\tg(L)=\frac12$, where
\be
\tg(L) = \sup_{1\leq N\leq L-1} \sup_{f\in L^2(\nu)}\,
\frac{\var(f)}{\tilde\cD(f,f)} \,.
\la{tgl}
\end{equation}
As in the proof of Proposition \ref{pro_iter} --
see (\ref{iter60}) -- we obtain
\be
\nu(f(\id - P)f) \leq \frac{L-2}{L-1} \,\tg(L-1) \tilde\cD(f,f)\,.
\la{iter68}
\end{equation}
where again $P$ is defined by (\ref{iter30}) with $\cF_k$ the $\si$-algebra
generated by $\a_k\in\{0,1\}$. The analysis of the spectrum of $P$
is much simpler now. Indeed any $f$ of the form (\ref{eq1}), with $\nu(f)=0$,
here must be of type
\be
f(\a) = \sum_{\ell}a_\ell\bar \a_\ell\,,\quad\quad \bar \a_\ell :=\a_\ell -
\frac{N}{L}\,.
\la{eq00}
\end{equation}
As in (\ref{eq2}) a
simple computation shows that
if $f$ is given by (\ref{eq00}) then $(\id-P)f=\frac{L-2}{L-1}f$.
It follows that any $f\in L^2(\nu)$ such that $\nu(f)=0$ satisfies
\be
\nu(f(\id - P)f)\geq \frac{L-2}{L-1}\nu(f^2)
\la{eq05}
\end{equation}
By (\ref{iter68}) we then obtain
\be
\tg(L)\leq \tg(L-1)\leq \cdots \leq \tg(2)\,\quad\quad L\geq 3\,.
\la{tgtg}
\end{equation}
On the other hand simple computations show that whenever $N=1$
one has $$\tilde D(f,f)=2\var_\nu(f)$$ for any $f\in L^2(\nu)$. In particular
$\tg(2)=1/2$ and $\tg(L)\geq 1/2$ for all $L\geq 3$. Together with
(\ref{tgtg}) this shows that $\tg(L)=1/2$.
A more detailed spectral analysis of this model
including all the eigenvalues was obtained in \cite{DSh}
by a different technique.




\newpage
\section{Quantum XXZ Hamiltonian}
Given $S\in\frac12\bbN$, $H\in\bbN$, consider the Hilbert space
$\frH = \otimes_{h=1}^H\bbC^{2S+1}$.
The spin-S XXZ chain on $[1,H]\cap\bbZ$
with {\em kink} boundary conditions
is defined by the operator
\begin{gather}
\cH^{(S)}=\sum_{h=1}^{H-1} \cH_{h,h+1}^{(S)}\,,\la{ham}\\
\cH_{h,h+1}^{(S)}
= S^2 -\Delta^{-1}
\left(S^1_{h}S^1_{h+1} + S^2_{h}S^2_{h+1}\right)
-S^3_{h}S^3_{h+1} +S\sqrt{1-\Delta^{-2}}\left(S^3_{h+1}-S^3_h\right) \,.
\nonumber
\end{gather}
Here $S^i_h$, $i=1,2,3$, are the spin--S operators (the
$2S+1$--dimensional irreducible representation of $SU(2)$) at every $h$,
and the constant $S^2$ has been added in order to have zero ground state
energy.  The parameter $\Delta\in(1,\infty)$ measures the anisotropy
along the third axis.  The kink boundary condition is obtained through
the telescopic sum
$S^3_H-S^3_1=\sum_{h=1}^{H-1}\left(S^3_{h+1}-S^3_h\right)$ and the
pre-factor $S\sqrt{1-\Delta^{-2}}$ is chosen in order to obtain
non-trivial ground states describing quantum domain walls (see
\cite{Alcarazetal}, \cite{S} and references therein).  \\ We choose the
basis of $\frH$ labelled by the $2S+1$ states of the third component of
the spin at each site and we write it in terms of configurations
$$
m=(m_1,\dots,m_H)\in\{-S,-S+1,\dots,S-1,S\}^H=: \cQ_S\,
$$
so that $\ket{m}=\otimes_{h=1}^H \ket{m_h}$ stands for the generic
basis vector in $\frH$. With these notations, and introducing the
stair-operators $S^{\pm}=S^1 \pm i S^2$, the action of $S^i$, $i=1,2,3$,
is given by
\begin{gather}
S^3_h \ket{m_h}=m_h\ket{m_h}\,,\quad
S^{\pm}_h\ket{m_h}=
c_{\pm}(S,m_h)\ket{m_h\pm 1}\,.\la{spins}\\
c_{\pm}(S,m_h) :=
\sqrt{(S\mp m_h)(S\pm m_h +1)}\nonumber
\end{gather}
The action of $\cH^{(S)}$ is explicated
by rewriting the pair-interaction terms as
\be
\cH_{h,h+1}^{(S)}=
S^2 -(2\Delta)^{-1}
\left(S^+_{h}S^-_{h+1} + S^-_{h}S^+_{h+1}\right)
-S^3_{h}S^3_{h+1} +S\sqrt{1-\Delta^{-2}}\left(S^3_{h+1}-S^3_h\right)
\la{ham2}
\end{equation}

\subsection{The spectral gap}
The Hamiltonian $\cH^{(S)}$ commutes with the total third component
of the spin
$$
S^3_{\rm tot}=\sum_{h=1}^H S^3_h\,.
$$
We shall divide the
space $\frH$ into sectors $\frH_{n}$, $n\in\{-SH,-SH+1,\dots,SH-1,SH\}$,
given by the eigenspaces corresponding to the eigenvalue $n$ of
$S^3_{\rm tot}$. It is known \cite{Alcarazetal} that for each $n$
there is a unique (up to multiplicative constants)
vector $\psi_n\in\frH_{n}$ such that $\cH^{(S)}\psi_n=0$, which is
given by
\begin{gather}
\psi_n = \sumtwo{m\in\cQ_S:}{\sum_h m_h=n}
\psi(m) \ket{m} \nonumber\\
\psi(m) =
\prod_h q^{h m_h} \,\sqrt{\binom{2S}{S+m_h}}\,.
\la{ground}
\end{gather}
Here $q\in(0,1)$ is the anisotropy parameter linked to $\Delta$
by the equation
\be
\Delta = \frac12(q+q^{-1})\,.
\la{qD}
\end{equation}
The ground states $\psi_n$ are interpreted as describing an interface
profile, \cite{Alcarazetal,B1}.  A fundamental question associated to
the stability of such ``quantum interfaces'' is the positivity of the
spectral gap \cite{BK,BKS}. The latter, denoted $\gap(\cH^{(S)})$, is
defined as the energy of the first excited state, i.e.\ the first
non-zero eigenvalue of the non-negative operator $\cH^{(S)}$.  Recently
this question was studied in great detail in the paper \cite{BKS} by
both analytical and numerical means. One of the main results of
\cite{BKS} is a proof of the fact that for every $S\in\frac12\bbN$,
$\gap(\cH^{(S)})$ is positive uniformly in $H$. Furthermore it was
conjectured on the basis of numerical analysis that $\gap(\cH^{(S)})$
should grow linearly with $S$. We prove the following bounds.

\begin{Th}\la{gaps}
For every $\Delta\in(1,\infty)$, there exists $\delta > 0$ such that
$$
\delta\,S\leq\gap(\cH^{(S)})\leq \delta^{-1} S
$$
for all $S\in\frac12\bbN$ and all $H\geq 2$.
\end{Th}
In order to prove Theorem \ref{gaps}
we shall establish the following unitary equivalence.
Let $L=2S$ and $N=SH+n$, recall the definition of
$L^2(\hat \Om,\hat\nu)$, of the subspace of horizontally
symmetric functions and of the generator $\widehat\cL$ introduced in
(\ref{generas}).
The measure $\hat \nu$ in
(\ref{hatnuo}) can be written, using (\ref{ground}) with $m=\o-S$, as:
\be
\hat\nu(\o) = \frac1{\tilde Z} \,\big[\psi(\o - S)\big]^2
\,\quad\quad \tilde Z = \sumtwo{\o\in\hat\Om:}{\sum_h\o_h = SH+n}
\big[\psi(\o - S)\big]^2
\,.
\la{hatnu}
\end{equation}
For any $\varphi\in \frH_n$ we also write
$$
\varphi = \sumtwo{m\in\cQ_S:}{\sum_hm_h=n}
\varphi(m)\ket{m}\,.
$$
Finally we set $\tfi(\o)=\varphi(\o-S)$.
Then the transformation
$$
\varphi(m) \:\to\: \frac{1}{\sqrt{\hat\nu(\o)}}\:\tfi(\o)\,=:\,
\big[U_n\varphi\big](\o)\,,
\quad\quad \o = m+S
$$ maps unitarily $\frH_n$ into $L^2(\hat \Om,\hat\nu)$.

\begin{Le}
\la{unitaryo}
For every $n\in\{-SH,-SH+1,\dots,SH-1,SH\}$
\be
U_n \,\cH^{(S)}\,\varphi = -\frac{S}{\Delta}\,\widehat\cL \,U_n\,\varphi\,,
\quad\varphi\in\frH_n
\la{unitaryo1}
\end{equation}
\end{Le}
\proof
From (\ref{ham2}) we compute
\begin{align*}
&\cH^{(S)}_{(h,h+1)}\ket{m} =
\big(S^2-m_hm_{h+1}+S\sqrt{1-\Delta^{-2}}(m_{h+1}-m_{h})\big)\ket{m}\\
& \quad -(2\Delta)^{-1}c_+(S,m_h)c_-(S,m_{h+1})|m^{+,h}\rangle
-(2\Delta)^{-1}c_-(S,m_h)c_+(S,m_{h+1})|m^{-,h}\rangle\,.
\end{align*}
Here we are using the notation
$$
m^{\pm,h}_{h'} = \begin{cases}
m_{h'} & h'\neq h,h+1 \\
m_{h}\pm 1 & h'=h \\
m_{h+1} \mp 1 & h'=h+1
\end{cases}
$$
Therefore
\begin{align*}
& (2\Delta)\big[\cH^{(S)}_{(h,h+1)}\varphi\big](m) =
\Gamma(m_h,m_{h+1})\,\varphi(m) \\
& \quad \quad- c_+(S,m_h)c_-(S,m_{h+1})
\varphi(m^{+,h}) - c_-(S,m_h)c_+(S,m_{h+1})\varphi(m^{-,h})
\end{align*}
with 
$$
\Gamma(m_h,m_{h+1}) =
(2\Delta)\big[S^2-m_hm_{h+1}+S\sqrt{1-\Delta^{-2}}(m_{h+1}-m_{h})\big]\,.
$$
Now a computation shows that
$$
\Gamma(m_h,m_{h+1}) = w_{+,h}(\o) + w_{-,h}(\o)\,,\quad\quad \o=S+m\,,
$$
with $w_{\pm,h}$ the rates defined in (\ref{generas}).
Another computation shows that
$$
c_+(S,\o_h-S)c_-(S,\o_{h+1}-S)
\sqrt{\frac{\hat\nu(\o^{+,h})}{\hat\nu(\o)}}
 = w_{+,h}(\o)
$$
and similarly
$$
c_-(S,\o_h-S)c_+(S,\o_{h+1}-S)
\sqrt{\frac{\hat\nu(\o^{-,h})}{\hat\nu(\o)}}
 = w_{-,h}(\o)\,.
$$
We have then obtained
\begin{align*}
(2\Delta)\big[U_n\,\cH^{(S)}_{(h,h+1)}\varphi\big](\o) &=
\, \big(w_{+,h}(\o) + w_{-,h}(\o)\big) \big[U_n\varphi\big](\o)\\
& - w_{+,h}(\o) \big[U_n\varphi\big](\o^{+,h})
- w_{-,h}(\o) \big[U_n\varphi\big](\o^{-,h})
\end{align*}
and the lemma follows.
\qed

\\
We are now able to finish the proof of Theorem \ref{gaps}.
Recall (\ref{dirD0}). Since here $L=2S$,
we readily infer from Lemma \ref{unitaryo} the estimate
\be
\gap(\cH^{(S)}) \geq \frac{S}{\Delta\gamma(2S,H)}\,,
\la{gaps1}
\end{equation}
where $\gamma(2S,H)$ is defined by (\ref{gammalh}).
The bound $\gamma\leq \delta^{-1}$
is the content of Theorem \ref{main}. On the other hand to
prove $\gap(\cH^{(S)})\leq \d^{-1}\,S$
we may use the following simple argument
which says that in each sector $\frH_n$ there are excited states
with energy bounded by $\d^{-1}\,S$.
Choose $f(\o) = \o_{h_0}$ with $h_0 = [\r]$, $\r=N/(2S)$.
A simple estimate shows that
$$
\hat \nu\big(f(-\widehat\cL)f\big) \leq
q^{-1}\,\hat\nu\big(\o_{h_0}+\o_{h_0+1}\big)
+ q\,\hat\nu\big(\o_{h_0}+\o_{h_0-1}\big) \leq k\,S\,(1\wedge\r)
$$
for some finite $k=k(q)<\infty$. On the other hand
using the estimates in \cite{CapMar} it is possible to check that
$\var_{\hat\nu}(f) \geq \d\, S\,(1\wedge\r)$ for some
$\d=\d(q)>0 $. By the variational principle
$$
\gap (\cH^{(S)}) \leq \,\frac{S}{\Delta}\,
\frac{\hat \nu\big(f(-\widehat\cL)f\big)}{\var_{\hat\nu}(f) }
\leq k \,S\,.
$$

\newpage
\section{Energy gap above the diagonal interfaces of the XXZ model}
In this section we study a higher dimensional
quantum $XXZ$ Hamiltonian which is
sometimes used to model a tilted interface (\cite{B1,BCNS,CapMar}).
In order to avoid complicated notation we shall work in a two-dimensional
setting. We later observe that our results actually hold
without modifications in any dimension.

\\
We begin by defining a cylindrical region (in two dimensions simply a rectangle)
with axis along the diagonal. Given two integers $R$ and $H$ we define
\be
\G = \G_{R,H} = \big\{x\in\bbZ^2:\; -R\leq x_1 - x_2 \leq R\,,\;
1\leq x_1+x_2\leq H \big\}
\la{grh}
\end{equation}
We write $\ell_x = x_1 + x_2$ for the distance of a site $x$
form the line $x_1 = -x_2$. A bond is an oriented pair
$b=(x,y)$, with $x,y\in\bbZ^2$ such that $|x_1-y_1|+|x_2-y_2|=1$
(in particular $\ell_y = \ell_x \pm 1$).
We call $\cB=\cB_{R,H}$ the set of bonds
$b=(x,y)$ with $x,y\in\G$ and $\ell_y = \ell_x + 1$.
For any $S\in\frac12\bbN$ the anisotropic
spin--$S$ Hamiltonian in the region $\G$ with kink boundary conditions
is defined by
\begin{gather}
\cH_{R}^{(S)}= \sum_{b\in\cB} \cH_b^{(S)}
\,,\la{rspins}\\
\cH_{b}^{(S)}
= S^2 -\Delta^{-1}
\left(S^1_xS^1_y + S^2_xS^2_y\right)
-S^3_xS^3_y
 + S\sqrt{1-\Delta^{-2}}\left(S^3_y-S^3_x\right) \,,\quad\; b=(x,y)
\nonumber
\end{gather}
with $S^k_x$, $k=1,2,3$, the spin--$S$ operators
at site $x$. This is the higher--spin analog of the spin--$\frac12$
cylindrical models considered in \cite{BCNS,CapMar}.
As usual we consider the Hilbert space
$\frH_\G = \otimes_{x\in\G}\bbC^{2S+1}$ with the basis labelled by configurations
$$
\ket{m} = \otimes_{x\in\G}\ket{m_x}\,\quad
m\in \{-S,\dots,S\}^\G=:\cQ_{\G,S}\,.
$$
Clearly the total third component
$S^3_{\rm tot}=\sum_{x\in\G} S^3_x$ is conserved and we may divide
$\frH_\G$ into sectors $\frH_{\G,n}$, $n\in\{-S|\G|,\dots,S|\G|\}$
according to the eigenvalues of $S^3_{\rm tot}$. Following \cite{Alcarazetal}
we know that in each such sector there is a unique ground state $\psi_{\G,n}$
given by
\begin{gather}
\psi_{\G,n} = \sumtwo{m\in\cQ_{\G,S}:}{\sum_x m_x=n}
\psi_{\G}(m) \ket{m} \nonumber\\
\psi_\G(m) =
\prod_{x\in\G} q^{\ell_x  m_x} \,\sqrt{\binom{2S}{S+m_x}}\,.
\la{grounds}
\end{gather}
As before $q=q(\Delta)$ is defined by (\ref{qD}).
The above ground states have zero energy
and we call $\gap(\cH_{R}^{(S)})$ the
first nonzero eigenvalue of $\cH_{R}^{(S)}$.

\\
Our main result here is a generalization to higher--spin
models of a theorem we proved in the case $S=\frac12$, \cite{CapMar}.
\begin{Th}
\la{hd}
For every $\Delta\in(1,\infty)$, there
exists $\delta > 0$ such that for all $S\in\frac12\bbN$, all $R,H\in\bbN$
\be
\delta\,S\,R^{-2}\leq\gap(\cH_{R}^{(S)})\leq \delta^{-1}\,S\,R^{-2}\,.
\la{hd1}
\end{equation}
\end{Th}

\begin{remark}
It will be clear from the proof that exactly the same estimates
hold (except that now the constant $\d$ in (\ref{hd1})
may depend on the dimension)
in the case where $\G$ is replaced by
a cylinder in $d+1$ dimensions with axis along the $11\cdots 1$
direction and basis given by a $d$-dimensional hypercube of side $R$.
\end{remark}

\proof
As a first step we establish a unitary equivalence
in the spirit of Lemma \ref{unitaryo}. Namely let
$\tilde\Om_\G = (2S+1)^\G$ denote the set of configurations
$\o_x\in\{0,\dots,2S\}$, $x\in\G$. Following our previous analysis
we may interpret $\o_x$ as the number of particles at site $x$.
Given $N\in\{0,1,\dots,2S|\G|\}$ define
$$
\tilde\Om_{\G,N} = \big\{
\o\in\tilde\Om_\G:\; \sum_{x\in\G}\o_x = N\,\big\}
$$
Consider
the probability measure on $\tilde\Om_{\G,N}$ given by
\be
\hat{\nu}_{\G,N}(\o) = \frac1{Z_{\G,N}}\big[\Psi_\G(\o - S)\big]^2\,,
\quad Z_{\G,N}=\sum_{\o\in\tilde\Om_{\G,N}} \big[\Psi_\G(\o - S)\big]^2\,.
\la{hatnugn}
\end{equation}
Then as in
Lemma \ref{unitaryo} we
obtain that
$$
\varphi(m) \:\to\: \frac{1}{\sqrt{\hat\nu_{\G,N}(\o)}}
\:\varphi(\o-S)\,=:\,
\big[U_{\G,n}\varphi\big](\o)\,,
\quad\quad \o = m+S
$$
maps unitarily each sector $\frH_{\G,n}$
into $L^2(\tilde\Om_{\G,N},\hat\nu_{\G,N})$, $N=S|\G|+n$.
Repeating the computation leading
to (\ref{unitaryo1}), for every bond $b\in\cB$ we have
\be
U_{\G,n}\,\cH^{(S)}_b\,U_{\G,n}^{-1}\,f
= -\,\frac1{2\Delta}\, \big\{w_b^+ \grad_b^+f + w_b^- \grad_b^-f\big\}\,,\quad
f\in L^2(\tilde\Om_{\G,N},\hat\nu_{\G,N})
\la{hd2}
\end{equation}
with the notation
\begin{gather*}
w_b^+(\o)= q^{-1}\,\o_y(2S-\o_x)\,,\;w_b^-(\o)= q\,\o_x(2S-\o_y)\,,\quad
b=(x,y)\,,\\
\grad_b^\pm f(\o) = f(\o^{b,\pm})-f(\o)\,,\quad \o^{b,\pm}_z :=
\begin{cases}
\o_z & z\neq x,y\\
\o_x \pm 1 & z=x \\
\o_y \mp 1 & z=y
\end{cases}
\end{gather*}
Let us define the Markov generator $\widehat\cG_{R,S}$
by
\be
\widehat\cG_{R,S}f(\o) = \sum_{b\in\cB}\big\{w_b^+(\o) \grad_b^+f(\o)
+ w_b^-(\o) \grad_b^-f(\o)\big\}\,.
\la{gengrs}
\end{equation}
Then $\widehat\cG_{R,S}$ is symmetric in
$L^2(\tilde\Om_{\G,N},\hat\nu_{\G,N})$ and by (\ref{hd2})
we have the unitary
equivalence
\be
\cH_R^{(S)} \simeq \,-\,\frac1{2\Delta}\,\widehat\cG_{R,S}
\la{hd4}
\end{equation}
Set $\Om_{\G^{(S)}} = \{0,1\}^{\G^{(S)}}$, where $\G^{(S)}$ is
the $3$--dimensional region
$$
\G^{(S)} = \big\{(i,x)\,:\quad i=1,\dots,2S\,;\;x\in\G\big\}\,.
$$
Define then $\Om_{\G^{(S)},N}$ as the set of $\a\in\Om_{\G^{(S)}}$ such that
$\sum_{(i,x)\in\G^{(S)}}\a_{(i,x)}=N$. 
Consider the probability measure $\nu_N$ on $\Om_{\G^{(S)},N}$
defined by
\be
\nu_N(\a) = \frac1{Z_N}
\prod_{(i,x)\in\G^{(S)}} q^{2\ell_x\a_{(i,x)}}\,,
\quad  Z_N = \sum_{\a\in\Om_{\G^{(S)},N}}
\prod_{(j,y)\in\G^{(S)}} q^{2\ell_y\a_{(j,y)}}
\la{hd5}
\end{equation}
Consider the subspace $\cS_{\G}\subset L^2(\Om_{\G^{(S)},N},\nu_N)$
of symmetric functions defined by: $f(\a)=f(\a^{\pi,x})$ for all $\a$,
all $x\in\G$ and all permutations $\pi\in\cP_{2S}$,
with $\a^{\pi,x}$ defined as in (\ref{symmpi}) replacing $h$ by $x$.
As in section \ref{exclusion},
$\hat{\nu}_{\G,N}$ can be looked at as the marginal
of $\nu_N$ on the sums
$$
\o_x = \sum_{i=1}^{2S}\a_{(i,x)}\,,\quad\quad
\a\in\Om_{\G^{(S)},N}
$$
and $\cS_\G$ is identified with $L^2(\tilde\Om_{\G,N},\hat\nu_{\G,N})$.
Then $\widehat\cG_{R,S}$ may be identified with the restriction to
$\cS_{\G}$ of $(2S)\,\cG_{R,S}$, the Markov generator defined by
\begin{gather}
\cG_{R,S}f(\a) = \frac1{2S}\,\sum_{i=1}^{2S}\sum_{j=1}^{2S}
\sum_{(x,y)\in\cB} c_{(i,x);(j,y)}(\a)\,\grad_{(i,x);(j,y)}f(\a)
\la{hd6}\\
c_{(i,x);(j,y)}(\a) = q^{\al_{(i,x)}-\al_{(j,y)}}\,,\quad
\grad_{(i,x);(j,y)}f(\a) = f(\al^{(i,x);(j,y)})-f(\al)\,.
\nonumber
\end{gather}
Here as usual $\al^{(i,x);(j,y)}$ denotes the configuration
after the exchange between $(i,x)$ and $(j,y)$. At this point we have
obtained a unitary equivalence
\be
\cH_R^{(S)} \simeq \,-\,\frac{S}{\Delta}\,\cG_{R,S}
\la{hd7}
\end{equation}
when the right hand side above is restricted to $\cS_\G$.
Notice the analogy of (\ref{hd6}) with the process introduced in
section \ref{exclusion}, see (\ref{genera}). However there is an important
difference (which will be seen in a moment to be responsible for
the $R^{-2}$ factors in (\ref{hd1})): while particles diffuse asymmetrically
in the $11$ direction just as it happens for (\ref{genera})
along the vertical
direction, we have in (\ref{hd6}) in addition
an essentially symmetric diffusion
along the orthogonal direction (given by the lines $\ell_x=$ constant).
We are now going to take care of these facts.

\smallskip

Observe that $\cG_{R,S}$ is symmetric in $L^2(\Om_{\G^{(S)}},\nu_N)$
with Dirichlet form
\begin{gather}
\cE_{R,S}(f) = \nu_N\big(f(-\cG_{R,S})f\big) =
 \frac1{2S}\,\sum_{i=1}^{2S}\sum_{j=1}^{2S}\, E_{(i,j)}(f)
\la{hd8}
\\
E_{(i,j)}(f):= \frac12 \sum_{(x,y)\in\cB}
\nu_N\,\Big[c_{(i,x);(j,y)}\big(\grad_{(i,x);(j,y)}f\big)^2
\Big]
\nonumber
\end{gather}
For any $x\in\G$ let $t_x=x_1-x_2$, so that $x$ is uniquely
determined by a couple $(t,\ell)$, $-R\leq t \leq R$,
$1\leq \ell\leq H$.
We may then write
\be
E_{(i,j)}(f) = \frac12
\sumtwo{t,s\in [-R,\dots,R]:}{|t-s|=1}\sum_{\ell=1}^{H-1}
\chi(t,\ell)\,
\nu_N\,\Big[c_{(i,t,\ell);(j,s,\ell+1)}\big(\grad_{(i,t,\ell);(j,s,\ell+1)}
f\big)^2
\Big]
\la{hd99}
\end{equation}
where $\chi(t,\ell)=1$ if there exists $x\in\G$ such that $x_1-x_2=t$
and $x_1+x_2=\ell$ and $\chi(t,\ell)=0$ otherwise
(notice that if $\chi(t,\ell)=1$ then
$\chi(s,\ell+1)=1$ for $s\in [-R,\dots,R]$ with $|s-t|=1$).
Using the bounds $q \leq c \leq q^{-1}$ on the rates
and the properties of
the measure $\nu_N$ it is not difficult
(see Lemma 2.6 in \cite{CapMar} for similar
computations) to show that there exists
$k=k(q)<\infty$ such that
given arbitrary $t,s\in [-R,\dots,R]$ with $t\leq s-1$ and
$\ell\in[1,H]$ one has
\begin{align}
\chi(t,\ell)&\chi(s,\ell+1)
\,\nu_N\,\Big[c_{(i,t,\ell);(j,s,\ell+1)}\big(\grad_{(i,t,\ell);(j,s,\ell+1)}
f\big)^2 \Big] \nonumber\\
& \leq \,k\,R \sum_{r=t}^{s-1}
\chi(r,\ell)\,\nu_N\,\Big[c_{(i,r,\ell);(j,r+1,\ell+1)}\big(\grad_{(i,r,\ell);(j,r+1,\ell+1)}
f\big)^2
\Big]
\la{hd10}
\end{align}
From (\ref{hd99}) and (\ref{hd10}) we easily obtain
\begin{align}
\sum_{t,s\in [-R,\dots,R]}
\sum_{\ell=1}^{H-1}
\chi(t,\ell)\chi(s,\ell+1)
&\,\nu_N\,\Big[c_{(i,t,\ell);(j,s,\ell+1)}\big(\grad_{(i,t,\ell);(j,s,\ell+1)}
f\big)^2
\Big] \nonumber\\
& \leq \,k\,R^3\,E_{(i,j)}(f)
\la{hd11}
\end{align}
Recalling (\ref{hd8}) we may summarize the above
estimate with the following statement
\begin{gather}
\cE_{R,S}(f)\geq \,\d\, R^{-2} \cD_{R,S}(f)\,,\la{hd100}\\
\cD_{R,S}(f) := \frac1{2SR}
\sum_{i,j\in[1,\dots,2S]}\sum_{t,s\in [-R,\dots,R]}
D_{(i,t);(j,s)}(f)\,,\nonumber\\
D_{(i,t);(j,s)}(f):= \sum_{\ell=1}^{H-1}
\chi(t,\ell)\chi(s,\ell+1)
\,
\nu_N\,\Big[c_{(i,t,\ell);(j,s,\ell+1)}\big(\grad_{(i,t,\ell);(j,s,\ell+1)}
f\big)^2 \Big]
\nonumber
\end{gather}
We are now able to conclude thanks to Theorem \ref{main}.
Indeed, the Dirichlet form $\cD_{R,S}$ defined above is a
special case of the one appearing in the Theorem, when we set
$L=2SR$ if $R$ is even and $L=2S(R+1)$ if $R$ is odd.
We then have that
$$
\cD_{R,S}(f) \geq \,\d \,\var_{\nu_N}(f)
$$
for some uniform constant $\d$. From (\ref{hd100}) and (\ref{hd7})
we see that the energy of excited states in each sector
lies above $\d \, S\, R^{-2}$ for some uniform constant $\d > 0$.
The desired lower bound on $\gap(\cH^{(S)}_R)$ follows at once.
The reverse estimate is much easier and it
can be obtained again from (\ref{hd7})
by modifying slightly the reasoning in Proposition 6.1 of \cite{CapMar}.

\qed

\newpage



%
%
%
%
%

\end{document}